\documentclass[12pt,reqno]{article}

\usepackage{amssymb,amsfonts,amsmath}
 
\usepackage[a4paper]{geometry}

\usepackage{parskip}
 
\usepackage{graphicx}
 
\usepackage{textcomp}
 
\newtheorem{theorem}{Theorem}[section]
\newtheorem{lemma}[theorem]{Lemma}

\usepackage{mathtools}
\DeclarePairedDelimiterX{\norm}[1]{\lVert}{\rVert}{#1}

\usepackage{fancyhdr}
\title{Rotation Sets of Open Billiards}
\author{Zainab Alsheekhhussain}
\begin{document}
\def\endofproof{{\rule{6pt}{6pt}}}

\titlepage
\maketitle

\begin{abstract}
We investigate the rotation sets of open billiards in $\mathbb{R}^N$ for the natural observable related to a starting point of a given billiard trajectory. We prove that the general rotation set is convex and the set of all convex combinations of rotation vectors of periodic trajectory $P_\phi$ is dense in it. We provide a constructive proof which illustrates that the set $P_\phi$ is  dense in the pointwise rotation set, and the closure of the pointwise rotation set is convex. We also  consider a class of billiards consisting of three obstacles and construct a sequence  in the symbol space such that its rotation vector is not defined.
\end{abstract}



\let\thefootnote\relax\footnotetext{Financial support by the Ministry of Higher Education in Saudi Arabia is greatfully acknowledged.}

\section{Introduction}
Billiards, traditionally, have been studied as part of Ergodic Theory, where their statistical properties have been investigated with respect to the invariant measure equivalent to the Lebesque measure. From that point of view, the limit behaviour of almost all trajectories has been studied, but it is also important to investigate the limit behaviour of all trajectories, specifically, \emph{periodic trajectories} (which are of zero measure). In this context, billiards in convex domains are considered twist maps, so elements of the so called Rotation Theory can be applied (see \cite{6}, \cite{24}, \cite{25}).

The dynamical model that we study here is the billiard flow in the exterior of several disjoint and strictly convex bodies with smooth boundaries satisfying a standard no eclipse condition (the convex hull of any two of the convex bodies has no common points with any of the other bodies \cite{3}). Such a system is called an \emph{open billiard}. This is quite a different system since the most interesting part of its phase space -- the so called non-wandering set -- is very small. It is locally homeomorphic to a Cantor set, and its Lebesgue measure is zero \cite{3}, \cite{8}. We associate with it an observable that is related to a starting point of a given billiard trajectory.

In general, Rotation Theory studies a given dynamical system and associates with it a function on the phase space with values in a vector space -- the so called \emph{observable}. Then by taking limits of ergodic averages of the observable, we get the \emph{rotation vectors (numbers)}. These limits (rotation vectors) form the \emph{general rotation set}. We can get also a rotation set by integrating the observable with respect to all ergodic invariant probability measures using Birkhoff's Ergodic Theorem \cite{11}, \cite{5}.

It is well-known that, topologically, an open billiard (considered as a discrete dynamical system) is isomorphic to a Markov shift of finite type. For such shifits a general theory of rotation sets has been developed by Ziemian in \cite{11}, where she studied the rotation set in the case when the dynamical system is a transitive subshift of finite type with an observable that depends only on cylinders of length two. Most recently Kucherenko and Wolf studied the rotation set of a compact metric space together with an $m$-dimensional continuous potential in \cite{5}. They considered the rotation set defined by means of integrals of Borel invariant probability measures. Their main result is that every compact and convex subset of $\mathbb{R}^m$ is attained as a rotation set of a particular set of potentials within a particular class of dynamical systems. However, as demonstrated by Blokh, Misiurewicz and Simany in \cite{12}, in physical systems such as billiards, for some naturally defined observables, one can obtain more significant information. In fact, \cite{12} studies one of the simplest possible cases, namely, a billiard on a torus in the exterior of one single small convex obstacle and for one particular observable -   the so called displacement. Even in this (simple at a first glance) case significant difficulties appear.

Let $\phi$ be the observable in the phase space of an open billiard related to a starting point in a billiard trajectory (see section 2 for details).
Let $P_\phi$ be the set consisting of all convex combinations of the form $conv(\rho_1,\rho_2,...,\rho_k)$, where $\rho_i=\rho_\phi(\xi_i)$ is the rotation vector of a periodic trajectory $\xi_i$. The folowing theorem summarizes our main result:
\begin{theorem}
\label{12}
The general rotation set for $\phi$ is convex, and  $P_\phi$  is dense in it.
\end{theorem}
This is derived from a result of Ziemian \cite{11} using some general arguments from the ergodic theory of sub-shifts of finite type. Moreover, the argument is fairly general and works for observables on any sub-shift of finite type.

Most of the paper deals with constructive arguments concerning the particular case of an open billiard. We first illustrate that the closure of the pointwise rotation set is convex. Then, in section 4, we consider a class of billiards consisting of three obstacles and construct a sequence $\xi$ in the symbol space such that the rotation vector $\rho_\phi(\xi)$ is not defined (see section 2). This can be easily generalized to  more general cases of open billiards.

The methodology used in this paper can be used to derive similar results for various other physical observables defined on an open billiard.
\section{Open Billiard}
Let $\Omega$ be a domain in $\mathbb{R}^N, N\geq 2$,  with smooth boundary $X=\partial\Omega$. The dynamical system generated by the motion of a point particle in $\Omega$ is called the \emph{billiard flow}. The point moves at constant velocity in the interior of the domain $\Omega$ creating reflections at the boundary $\partial\Omega$ according to the classical law of geometrical optics: the angle of incidence is equal to the angle of reflection (\cite{3}, \cite{7}, \cite{8}).

The open billiard is a kind of billiard in which $\Omega$ is unbounded and $\Omega=\mathbb{R}^N\backslash K$, where K is a union $K=K_1\cup K_2\cup ....\cup K_s$ of pairwise disjoint compact and strictly convex sets with $C^2$-smooth boundaries $\partial K_i$, for some $s\geq 3$. The obstacles $K_i$'s satisfy the non-eclipse condition. That is, for any $i\neq k\neq j$, the convex hull of $K_i\cup K_j$ does not intersect $K_k$.
The \emph{billiard ball map}, $B$ is defined on a subset $M'$ of the set 
$$M=\{(x,v)\in X\times {\bf{S}}^{N-1}: \langle\nu(x),v\rangle>0\},$$
where $\nu(x)$ is the unit normal to $\partial\Omega$ at $x$ pointing into the interior of $\Omega$. Let $q=(x,v)$ denote the trajectory $\gamma$ which starts at $x$ with direction $v$ and has a common point with $X$. Let $y$ denote the first such point. This means that $y\in X$ and the open segment $(x,y)$ has no common points with $X$. Define $M'$ to be the subset of $M$ consisting of all $q=(x,v)$ such that $\gamma$ issued from $x$ in direction $v$ hits $\partial \Omega$ at, say, $y$. Define $u$ to be $u=v-2\langle v,\nu(y)\rangle\nu(y)$ and set $(x,u)\xrightarrow{B}(y,v)$. So we obtain a map $B:M'\longrightarrow M$. Here we consider $B:M_0\longrightarrow M_0$, where $M_0$ is the set of all $(x,v)\in M'$ such that $B^n(x,v)$ is defined for all $n\in\mathbb{Z}$, which is called the \emph{non-wandering set}. We call the points $q\in M_0$ with $B^k(q)=q$ for some $k>0$ \emph{periodic points} of period $k$ of $B$.

The set $M_0$ can be coded using symbolic dynamics when the non-eclipsing condition holds. This is done by using the symbol space $\Sigma$ of infinite addmissible sequences $\xi=(....,\xi_{-1},\xi_0,\xi_1,....)$, where each $\xi_i\in \{1,2,...,s\}$. An $\emph{admissible word}$ (or sequence) is, by definition, a finite or an infinite sequence $\{\xi_i\}_{i\in I}$, where $I$ is an interval in $\mathbb{Z}$, i.e. $I=[p,q]$, for some $p<q$, such that $\xi_i\neq\xi_{i+1}$ for all $p\leq i<q$ (this includes the cases when $p=-\infty$ and/or $q=\infty$). These sequences are acted upon by the left-shift map $\sigma:\Sigma\longrightarrow\Sigma$, given by $\sigma(\xi)=\xi'$, where $\xi'=(\xi'_i)$ is defined by $\xi_i'=\xi_{i+1}$ (\cite{29}).

The observable $\phi$, considered here, is related to a starting point of a given billiard trajectory defined as follows: given $\xi=(\xi_i)$ in the symbol space $\Sigma$, there exists a unique pair $(y,v)$ of a point $y$ on the boundary of $K$ and a unit vector $v$ in $\mathbb{R}^N$ such that the billiard trajectory issued from $y$ in direction $v$ has infinitely many reflections both backwards and forwards and the i'th reflection occurs at the boundary of the obstacle $K_{\xi_i}$ for every integer $i$ (see \cite{9}, \cite{3}, \cite{8}). Then define $\Phi:\Sigma\longrightarrow M_0$ by $\Phi(\xi)=(y,v)$ and $\phi:\Sigma\longrightarrow\mathbb{R}^N$ by $\phi(\xi)=y$. Then $B(\Phi(\xi))=\Phi(\sigma(\xi))$ for all $\xi\in\Sigma$. The set of all limits 
$$\rho_\phi(\xi)=\lim_{n\rightarrow\infty}\frac{1}{n}\sum_{i=0}^{n-1}\phi(\sigma^i(\xi))$$ 
(for those $\xi\in\Sigma$ for which the limit exists) is called the \emph{pointwise rotation set} and will be denoted by $J_\phi$. The set $G_\phi$ of all limits of subsequences of sequences as in the right-hand-side of the above formula for all $\xi\in\Sigma$ is called the \emph{general rotation set}.

We will assume that the set $K$ is contained in a ball with radius $R$, so that $\norm{x}\leq R$ for every $x\in K$.

\section{The General Rotation Set}
\setcounter{equation}{0}
\renewcommand{\theequation}{\arabic{section}.\arabic{equation}}
In this section, we will show that the general rotation set $G_\phi$ is convex and equal to the closure of $P_\phi$. We will use results from \cite{11} to do that, but first we will need to approximate our observable $\phi$ with ones that are locally constant.

Note that $\phi$  is continuous (in fact, H\"{o}lder continuous) see Lemma \ref{5} below. 

Thus, we can use \emph{Stone-Weierstrass Theorem} to find approximations of $\phi$.  Let $C(\Sigma)$ be the set of all real valued continuous functions defined on $\Sigma$ and $$A=\{\varphi\in C(\Sigma)\,:\,\varphi\, depends\,  on\,  finitely\,  many\,  coordinates \},$$ then every constant function belongs to $A$, and if $\varphi$ and $\chi$ are in $A$ then $\varphi +\chi$  and $\varphi\chi$ are in $A$. Also for every $\varphi$ in $A$ and every $\lambda$ in $\mathbb{R}$, $\lambda\varphi$ is in $A$. 
Thus, $A$ is an algebra in $C(\Sigma)$ containing the constants. Also if $\xi,\xi'$ are elements in $\Sigma$ such that $\xi\neq\xi'$, then there exists $i$ such that $\xi_i\neq\xi_i'$. 
Consider $\varphi:\Sigma\longrightarrow\Sigma$ such that $\varphi(\eta)=\eta_i$ for every $\eta$ in $\Sigma$. (Here we assume that each $\eta_i=1,2,...,s$, so the values of $\eta_i$ are real numbers.) Then $\varphi$ is in $A$ since $\varphi$ depends on only one coordinate. Clearly $\varphi(\xi)\neq\varphi(\xi')$. 
Thus, $A$ is dense in $C(\Sigma)$ (i.e for every $\epsilon>0$ and every $\chi$ in $C(\Sigma)$, there exists $\varphi$ in $A$ such that $\parallel\chi-\varphi\parallel\leq\epsilon$).  Here $\norm{\varphi}=max_{\xi\in\Sigma}\mid\varphi(\xi)\mid$.

Next, we use a procedure described in \cite{28} (see section 10 in \cite{28} and Chapter 1 of \cite{30}).
 
Let $m\ge 2$. Let $F$  be the set of all blocks of the form $(\ell_1,\ell_2,...,\ell_m)=L$ such that $\ell_i\neq\ell_{i+1}$ and $\ell_i=1,2,...,s$. Now $F$ will be our new set of symbols. Define a new symbol space $\Sigma_F=\{ L=(L_i)_{i=-\infty}^\infty: L_i\in F \}$. Define $\Phi: \Sigma\longrightarrow\Sigma_F$ by $\Phi(\xi)=L=(L_i)_{i=-\infty}^\infty$, where $L_i=(\xi_i,\xi_{i+1},...,\xi_{i+m-1})$.  Set $\tilde{\Sigma}=\Phi(\Sigma)$. Hence, if $\varphi$ is in $C(\Sigma)$ then $\varphi=\psi\circ\Phi$ for some $\psi$ in $C(\tilde{\Sigma})$.

Now we want to show that $G_\varphi=G_\psi$. To do this, we need first to prove that $\Phi:\Sigma\longrightarrow\tilde{\Sigma}$ is bijective and $\Phi\circ\sigma^i=\sigma^i\circ\Phi$. 

Let $\xi,\tilde{\xi}$ be in $\Sigma$ such that $\xi\neq \tilde{\xi}$. (I.e $\xi_i\neq\tilde{\xi}_i$ for some $i$.) Assume that $\Phi(\xi)=L$ and $\Phi(\tilde{\xi})=\tilde{L}$. Then $L_i=(\xi_i,\xi_{i+1},...,\xi_{i+m-1})\neq \tilde{L}_i=(\tilde{\xi}_i,\tilde{\xi}_{i+1},....,\tilde{\xi}_{i+m-1})$. I.e $L\neq\tilde{L}$, so $\Phi$ is injective and, hence, $\Phi$ is bijective. 

Let $\xi=(.....,\xi_{-1};\xi_0,\xi_1,....)$ be in $\Sigma$. Then
 \begin{eqnarray}\Phi(\sigma^i(\xi))&=&\Phi((.....,\xi_{i-1};\xi_i,\xi_{i+1},....))\nonumber\\
&=&(.....,(\xi_{i-1},....,\xi_{i+m-2});(\xi_i,....,\xi_{i+m-1}),(\xi_{i+1},....,\xi_{i+m}),....)\nonumber\\
&=&\sigma^i((.....,(\xi_{-1},....,\xi_{m-2});(\xi_0,....,\xi_{m-1}),(\xi_{1},....,\xi_{m}),....))\nonumber\\
&=&\sigma^i(\Phi(\xi)).\nonumber
\end{eqnarray}

Now let $\rho$ be in $G_\varphi$. This means that there exists $\xi$ in $\Sigma$ such that $$\lim_{k\rightarrow\infty}\frac{1}{n_k}\Sigma_{i=0}^{n_k-1}\varphi(\sigma^i(\xi))=\rho$$ for some sequence $1\leq n_1\leq n_2\leq....$ .  Assume that $\Phi(\xi)=L$. So, for every $i$, $$\varphi(\sigma^i(\xi))=\psi\circ\Phi(\sigma^i(\xi))=\psi(\sigma^i(\Phi(\xi)))=\psi(\sigma^i(L)).$$ So $\rho$ is in $G_\psi$. I.e $G_\varphi$ is subset of $G_\psi$. Similarly, one can prove that $G_\psi$ is subset of $G_\varphi$. Hence, $G_\varphi=G_\psi$.

Now assume $\varphi$ depends on $m+1$ coordinates. Then $\psi$ depends on just two coordinates in $\tilde{\Sigma}$. Applying \cite{11} to $\psi$, we obtain a corresponding result about the general rotation set of $\varphi$. This follows from  $G_\varphi=G_\psi$. Thus, if $P_\psi$ is the set of all convex combinations of rotation vectors of periodic trajectories in $\tilde{\Sigma}$ then, by Ziemian \cite{11}, $P_\psi$ is dense in $G_\psi$, and, hence, $P_\varphi$ is dense in $G_\varphi$.

For our observable $\phi$ and $q>1$, using the above argument, we can find an observable $\psi=\psi_q$ depending only on finitely many coordinates such that $\norm{\phi-\psi}\leq\frac{1}{q}$. 

Now fix $\epsilon>0$ small, and fix $q>1$ a large integer such that \begin{equation}\frac{4}{q}<\epsilon\end{equation}
Let $\rho$ be in $G_\phi$, so $\rho$ has the form $\rho=\lim_{k\rightarrow\infty}\rho_k$ where $\rho_k=\frac{1}{n_k}\Sigma_{i=0}^{n_k-1}\phi(\sigma^i(\xi))$,  for some sequence of integers $1\leq n_1< n_2<....<n_k<...$ and some $\xi$ in $\Sigma$.  Set $\rho_k'=\frac{1}{n_k}\Sigma_{i=0}^{n_k-1}\psi(\sigma^i(\xi))$.
Then for any integer $k$:
\begin{eqnarray}
\parallel\rho_k-\rho_k'\parallel&=&\norm[\bigg]{\frac{1}{n_k}\Sigma_{i=0}^{n_k-1}\phi(\sigma^i(\xi))-\frac{1}{n_k}\Sigma_{i=0}^{n_k-1}\psi(\sigma^i(\xi))}\nonumber\\
&\leq&\frac{1}{n_k}\Sigma_{i=0}^{n_k-1}\norm{\phi(\sigma^i(\xi))-\psi(\sigma^i(\xi))}\nonumber\\
&\leq&\frac{1}{n_k}\Sigma_{i=0}^{n_k-1}\frac{1}{q}=\frac{1}{q}.\end{eqnarray}

Hence $\{\rho'_k\}$ is bounded, so it has a convergent subsequence, say, $\{\rho'_{k_m}\}$ going towards $\rho'$ in $G_\psi$ as $m$ goes to infinity. For every $m$, we have that $\rho'_{k_m}=\frac{1}{n_{k_m}}\Sigma_{i=0}^{n_{k_m}-1}\psi(\sigma^i(\xi))$. Set $\rho_{k_m}=\frac{1}{n_{k_m}}\Sigma_{i=0}^{n_{k_m}-1}\phi(\sigma^i(\xi))$. As $\{\rho_{k_m}\}$ is a subsequence of $\{\rho_k\}$, we have that $\{\rho_{k_m}\}$ converges to $\rho$. We also have by (3.2) that $\norm{\rho_{k_m}-\rho'_{k_m}}\leq\frac{1}{q}$. Thus, as m goes to infinity, we obtain \begin{equation}\norm{\rho-\rho'}\leq\frac{1}{q}\end{equation}
Similarly, we can prove that for any $\rho'$ in $G_\psi$, there exists $\rho$ in $G_\phi$ such that $\norm{\rho'-\rho}\leq\frac{1}{q}$. Hence the Hausdorff distance between $G_\phi$ and $G_\psi$ is less than or equal to $\frac{1}{q}$.

Let $\xi$ be a periodic trajectory in $\Sigma$ of period $p$. Then, by Lemma \ref{1} in section 4 below and the Remark, $\rho_\phi(\xi)=\frac{1}{p}\Sigma_{i=0}^{p-1}\phi(\sigma^i(\xi))$ and $\rho_\psi(\xi)=\frac{1}{p}\Sigma_{i=0}^{p-1}\psi(\sigma^i(\xi))$. Thus,
\begin{eqnarray}
\norm{\rho_\phi(\xi)-\rho_\psi(\xi)}&\leq&\frac{1}{p}\Sigma_{i=0}^{p-1}\norm{\phi(\sigma^i(\xi))-\psi(\sigma^i(\xi))}\nonumber\leq\frac{1}{q}.\end{eqnarray}

Let $P_\phi$ be the set consisting of all convex combinations of the form $conv(\rho_1,\rho_2,...,\rho_k)$, where $\rho_i=\rho_\phi(\xi_i)$ is the rotation vector of a periodic trajectory $\xi_i$. Let $P_\psi$  be the corresponding set for the observable $\psi$. If $\rho$ is in $G_\phi$ then there exists $\rho'$ in $G_\psi$ such that $\norm{\rho-\rho'}\leq\frac{1}{q}$, by (3.3). Since $\psi$ is locally constant, we can apply Theorem 3.4 of \cite{11} (see how above). Hence, $G_\psi$ is convex, and the set $P_\psi$  is dense in $G_\psi$. Thus $\rho'$ can be approximated by a  vector in $P_\psi$. I.e we can find $k$ positive real numbers $t_i$'s with $\Sigma_{i=1}^kt_i=1$ and rotation vectors $\rho_1',.....,\rho_k'$ of periodic trajectories  such that $\norm{\Sigma_{i=1}^kt_i\rho_i'-\rho'}\leq\frac{\epsilon}{2}$. Then $\norm{\Sigma_{i=1}^kt_i(\rho_i'-\rho_i)}\leq\Sigma_{i=1}^k\frac{t_i}{q}=\frac{1}{q}$. Clearly, $\Sigma t_i\rho_i$ is in $P_\phi$ and \begin{eqnarray}
\norm{\Sigma_{i=1}^kt_i\rho_i-\rho}&=&\norm{\Sigma_{i=1}^kt_i(\rho_i-\rho_i')+\Sigma_{i=1}^kt_i\rho_i'-\rho'+\rho'-\rho}\nonumber\\
&\leq&\frac{1}{q}+\frac{\epsilon}{2}+\frac{1}{q}\leq\epsilon.\nonumber\end{eqnarray}

Thus $P_\phi$ is dense in $G_\phi$.

It remains to prove that $G_\phi$ is convex.

Take $k$ elements,  $\rho_1,\rho_2,....,\rho_k$ of $G_\phi$. Set $\rho=\sum_{i=1}^kt_i\rho_i$, for some $t_1,t_2,...,t_k\geq 0$ and $\sum_{i=1}^kt_i=1$. Consider an arbitrary $\epsilon>0$ and take $q$ as in (3.1).  For each $\rho_i$, there exists $\rho_i'$ in $G_\psi$ such that $\norm{\rho_i-\rho_i'}\leq\frac{1}{q}$ by (3.3). As $G_\psi$ is convex, the sum $\Sigma_{i=1}^kt_i\rho_i'=\rho'$ is in $G_\psi$. Also, by (3.3) and the statement underneath, there is $\rho''$ in $G_\phi$ such that $\norm{\rho''-\rho'}\leq\frac{1}{q}$. So $$\norm{\rho-\rho''}\leq\norm{\rho-\rho'}+\norm{\rho'-\rho''}\leq\frac{2}{q}.$$ So, for every $q>1$, there exists $\rho_q''$ in $G_\phi$ such that $\norm{\rho-\rho_q''}\leq\frac{1}{q}$. Since $G_\phi$ is compact, $\{\rho_q''\}$ has a convergent subsequence, say $\{\rho_{q_i}''\}$ going towards $\rho'''$ (in $G_\phi$). From $\norm{\rho-\rho_{q_i}''}\leq\frac{2}{q_i}$, and as $i$ tends to infinity, we get $\norm{\rho-\rho'''}=0$. Thus, $\rho$ is in $G_\phi$, and hence, $G_\phi$ is convex.

This proves Theorem \ref{12}. \endofproof

\section{Periodic Orbits and the Pointwise Rotation Set}
\setcounter{equation}{0}
\renewcommand{\theequation}{\arabic{section}.\arabic{equation}}
In this section, we provide a constructive proof of the following weaker version of Theorem \ref{12}:
\begin{theorem}
\label{11}
On an open billiard, let $\phi$ be the observable defined by the position of a starting point of a given billiard trajectory in the symbol space. Then the closure of the pointwise rotation set for $\phi$ is convex, and the set $P_\phi$ is dense in the pointwise rotation set.
\end{theorem}
We need several lemmas in order to do that.

The periodic trajectories are generated by points $q\in M_0$ with $B^k(q)=q$ for some $k>0$. It is easy to see that all periodic trajectories have rotation vectors. 
\begin{lemma}
\label{1}
Let $\xi \in \Sigma$ be a periodic trajectory with period length $p$. Then $\xi$ has a rotation vector $\rho_\phi(\xi)=\frac{1}{p} \sum_{i=0}^{p-1}{x_i}$, where $x_i$ is the reflection point at the obstacle $K_{\xi_i}$. 
\end{lemma}
\emph{Proof}. We want to show that $\lim_{n\rightarrow\infty} \frac{1}{n}\sum_{i=0}^{n-1}\phi(\sigma^i\xi)$ exists. 

Let $\xi=(\xi_i)$ be periodic with period $p$, so $\xi_{i+p}=\xi_i$ for all $i$. Let $\phi(\sigma^i\xi)=(x_i,u_i)$, $(x_i,u_i)=(x_{i+np},u_{i+np})$ for all $n\in\mathbb{Z}$. We have $\phi (\sigma^i\xi) =x_i$. 

Let $n=kp+m,\, 0\le m\le p-1$, then
$$\sum_{i=0}^{n-1}\phi(\sigma^i\xi)=(k+1)(x_0+x_1+....+x_{m-1})+k(x_m+....+x_{p-1})$$
and so 
$$
\lim_{n\rightarrow\infty}\frac{1}{n}\sum_{i=0}^{n-1}\phi(\sigma^i\xi)=\lim_{k\rightarrow\infty}\frac{(k+1)(x_0+x_1+...+x_{m-1})+k(x_m+....+x_{p-1})}{kp+m}$$
$$=\frac{1}{p}\sum_{i=o}^{p-1}x_i$$
Thus $$\rho_\phi(\xi)=\frac{1}{p}\sum_{i=0}^{p-1}x_i.$$
And so every periodic trajectory has a rotation vector.
\endofproof 

\emph{Remark}:
The above lemma can be generalized for any observable.

\setlength{\parskip}{10pt plus 1pt minus 1pt}
To study rotation vectors of open billiards, we will need the following lemma which is a consequence of the strong hyperbolicity properties of the billiard flow - see \cite{10*}; see also Lemma 10.2.1 in \cite{8}.
\begin{lemma}
\label{5}
There exist constants $C>0$ and $\delta\in (0,1)$, depending only on $K$, with the following property: 

if $i_0,i_1,....,i_{m+1}\in\{1,2,...,s\}$ are such that $i_j\neq i_{j+1}$ for $0\leq j\leq m$,
$$x_0\in\partial K_{i_0}, x_1\in\partial K_{i_1},....., x_m\in\partial K_{i_m}, x_{m+1}\in\partial K_{i_{m+1}}$$
and
$$x_0'\in\partial K_{i_0}, x_1'\in\partial K_{i_1},.....,x_m'\in\partial K_{i_m}, x_{m+1}'\in\partial K_{i_{m+1}}$$
are two sequences of points such that for every $j=1,....,m$ the segments $[x_{j-1},x_j]$ and $[x_j,x_{j+1}]$ satisfy the law of reflection at $x_j$ with respect to $\partial K_{i_j}$ and the segments $[x_{j-1}',x_j']$ and $[x_j',x_{j+1}']$ satisfy the law of reflection at $x_j'$ with respect to $\partial K_{i_j}$, then
$$\parallel x_i-x_i'\parallel\le C(\delta^i+\delta^{m-i}),$$ for all $i=1,....,m$.
\end{lemma}

We will need the following set-up to prove the next three lemmas. Fix a small constant $\epsilon>0$ - we will say later how small $\epsilon$ should be. Then fix integers $\ell, p$ sufficiently large so that
\begin{equation}
\frac{5R}{p}<\frac{\epsilon}{3}, \; \delta^\ell<\frac{\epsilon}{3},
\end{equation}
where $\delta$ is a global constant in $(0,1)$ satisfying the conclusion of Lemma \ref{5}.

In the next lemma, we prove that every rotation vector of a periodic trajectory can be approximated by a rotation vector of another periodic trajectory of a special kind.

\begin{lemma}
\label{2}
For every periodic $\xi = (......;\overbrace{\xi_0,\xi_1,....,\xi_{n-1}}^{\mbox{$\Xi$}},.....)\in\Sigma$ of period $n$, there is a periodic $\tilde{\xi}$ in $\Sigma$ of the form $\tilde{\xi}=(......;\overbrace{\Xi,...,\Xi,....,\Xi}^{\mbox{$p\ell n$ times}},j,....)$, where $j\notin \{\xi_{n-1},\xi_0\}$ and $p$, $\ell$ are as in (4.1), such that $\parallel \rho_\phi(\xi)-\rho_\phi(\tilde{\xi})\parallel<\epsilon$.
\end{lemma}
\emph{Proof}. Let $\rho_\phi(\sigma^j(\xi))=x_j$, then, by Lemma \ref{1}, $\rho_\phi(\xi)=\frac{1}{n}\sum_{j=0}^{n-1}x_j$.
$\tilde{\xi}$ is also periodic of period $p\ell n+1$, and if $\rho_\phi(\sigma^j(\tilde{\xi}))=\tilde{x_j}$, then, by Lemma \ref{1}, $\rho_\phi(\tilde{\xi})=\frac{1}{p\ell n+1}\sum_{j=0}^{p\ell n}\tilde{x_j}$.
Clearly, in $\tilde{\xi}$ we can find at least $\ell n$ reflections from the same obstacles as in $\xi$ at the begining of the periodic block and at the end.
So, by Lemma \ref{5}, we have:
$$\parallel\tilde{x}_{\ell n+j}-x_j\parallel<\delta^{\ell n}+\delta^{\ell n}=2\delta^{\ell n},\; 0\le j\le(p-2)\ell n-1\;.$$
Thus,
\begin{eqnarray}
&&\norm[\bigg]{\sum_{j=\ell n}^{(p-1)\ell n-1}\frac{\tilde{x}_j}{p\ell n+1}-\sum_{j=0}^{n-1}\frac{x_j}{n}}\nonumber\\
&=&\norm[\bigg]{\frac{\sum_{i=0}^{(p-2)\ell-1}\sum_{j=0}^{n-1}\tilde{x}_{(\ell+i)n+j}}{p\ell n+1}-\frac{1}{n}\sum_{j=0}^{n-1}x_j}\nonumber
\end{eqnarray}
\begin{eqnarray}
&=&\norm[\bigg]{\frac{\sum_{i=0}^{(p-2)\ell n-1}\sum_{j=0}^{n-1}(\tilde{x}_{(\ell+i)n+j}-x_j)}{p\ell n+1}+\frac{(p-2)\ell}{p\ell n+1}\sum_{j=0}^{n-1}x_j-\frac{1}{n}\sum_{j=0}^{n-1}x_j}\nonumber\\
&\le&\frac{2\delta^{\ell n}(p-2)\ell n}{p\ell n+1}+\left|\frac{(p-2)\ell}{p\ell n+1}-\frac{1}{n}\right|\sum_{j=0}^{n-1}\parallel x_j\parallel\nonumber\\
&\le& 2\delta^{\ell n}+\frac{3R}{p}.
\end{eqnarray}
By (4.1) and (4.2) we have:
\begin{eqnarray}
&&\norm[\bigg]{\rho_\phi(\tilde{\xi})-\rho_\phi(\xi)}\nonumber\\
&=&\norm[\bigg]{\sum_{j=0}^{p\ell n}\frac{\tilde{x}_j}{p\ell n+1}-\sum_{j=0}^{n-1}\frac{x_j}{n}}\nonumber\\
&\le& \sum_{j=0}^{\ell n-1}\norm[\bigg]{\frac{\tilde{x}_j}{p\ell n+1}}+\norm[\bigg]{\sum_{j=\ell n}^{(p-1)\ell n-1}\frac{\tilde{x}_j}{p\ell n+1}-\sum_{j=0}^{n-1}\frac{x_j}{n}}+\sum_{j=(p-1)\ell n}^{p\ell n}\norm[\bigg]{\frac{\tilde{x}_j}{p\ell n+1}}\nonumber\\
&\le&\frac{2\ell n}{p\ell n+1}R+2\delta^{\ell n}+\frac{3R}{p}\nonumber\\
&\le&\frac{2R}{p+1/\ell n}+2\delta^\ell+\frac{3R}{p}<\frac{5R}{p}+2\delta^\ell<\frac{\epsilon}{3}+2\frac{\epsilon}{3}=\epsilon.\nonumber
\end{eqnarray}
This proves the lemma. \endofproof

Given points $\rho_1,\rho_2,....,\rho_k$ in $\mathbb{R}^N$, we will denote by $conv(\rho_1,\rho_2,....,\rho_k)$ the \emph{convex hull} of these points. i.e. the points $\upsilon$ of the form $\upsilon=\sum_{i=1}^kt_i\rho_i$ for some $t_1,....t_k\geq 0$ and $\sum_{i=1}^kt_i=1$.

The following lemma is the central point in this article.
\begin{lemma}
\label{3}
Let $\xi^{(1)}, \xi^{(2)},.....,\xi^{(k)}$ be periodic elements of $\Sigma$ and let $\rho_{i}=\rho_\phi(\xi^{(i)}), \, i=1,2,....,k$. Then for every $v\in conv(\rho_1,...,\rho_k)$ and every $\epsilon >0$, there exists a periodic $\eta\in \Sigma$ such that 
$\parallel\rho_\phi(\eta)-v\parallel<\epsilon$.
\end{lemma}
\emph{Proof}. Let $v=\sum_{i=1}^k t_i\rho_i$ for some $t_1,t_2,....,t_k >0$ and $\sum_{i=0}^k t_i=1$. We can approimate each $t_i$ by a rational number. There are positive integers $s_1,s_2,....,s_k$ such that \begin{equation}\norm[\bigg]{\frac{s_i}{m}-t_i}\le\frac{\epsilon}{3kR},\end{equation} where $m=\sum_{i=1}^k s_i$.

Let $\xi^{(i)}$ be as in Lemma \ref{2}. We construct $\tilde{\xi}^{(i)}$ by repeating $ps_i\ell_in_i$ times the periodic block of $\xi^{(i)}$, where $\ell_i=\frac{L}{n_i}$, $n_i$ is the length of the periodic block of $\xi^{(i)}$ and $L=\ell n_1n_2...n_k$, and $p$ and $\ell$ are as in (4.1). Applying Lemma \ref{2} to get $\tilde{\xi}^{(i)}$, define the last symbol $j_i=\tilde{\xi}^{(i)}_{m_i-1}\neq\tilde{\xi}^{i+1}_0$ for $i<k$ and $j_k\ne\tilde{\xi}^{(1)}_0$, where $$m_i =ps_i\ell_in_i=ps_iL.$$ By definition, $\tilde{\xi}^{(i)}$ is the periodic sequence obtained by repeating infinitely many times the block $\Xi^{(i)}=(\tilde{\xi}_0^{(i)},\tilde{\xi}_1^{(i)},....,\tilde{\xi}_{m_i-1}^{(i)})$ Then $$\tilde{\xi}^{(i)}=(....;\overbrace{\tilde{\xi}_0^{(i)},\tilde{\xi}_1^{(i)},....,\tilde{\xi}_{m_i-1}^{(i)}}^{\mbox{periodic bolck, $\Xi^{(i)}$}},....).$$
Assuming that $y_j^{(i)}, j=0,....,m_i-1$, are the reflection points of $\tilde{\xi}^{(i)}$ and $\phi(\sigma^j(\tilde{\xi}^{(i)}))=y_j^{(i)}$, then, by Lemma \ref{1}, as $\tilde{\xi}^{(i)}$ is periodic of length $m_i$, we have $$\rho_\phi(\tilde{\xi}^{(i)})=\frac{1}{m_i}\sum_{j=0}^{m_i-1}y_j^{(i)},$$ 
and \begin{equation}\parallel \rho_i-\rho_\phi(\tilde{\xi}^{(i)})\parallel<\frac{\epsilon}{3}.\end{equation}
Construct $\eta$ to be a periodic with periodic block as follows: $$\eta=(....;\overbrace{\Xi^{(1)},\Xi^{(2)},.............................,\Xi^{(k)}}^{\mbox{periodic block of length $\sum_{i=1}^km_i$}},....).$$ Let $z_j^{(i)},\; j=0,...,m_i-1, \; i=1,...,k$ be the reflection points of the corresponding billiards trajectory. Then $\eta$ is admissible, so $\eta\in\Sigma$.

Applying Lemma \ref{1}, we get $$\rho_\phi(\eta)=\frac{1}{\sum_{i=1}^km_i} \sum_{i=1}^k\sum_{j=0}^{m_i-1} z_j^{(i)}=\frac{1}{\sum_{i=1}^kps_iL} \sum_{i=1}^k \sum_{j=0}^{m_i-1}z_j^{(i)}=\frac{1}{mpL}\sum_{i=1}^k\sum_{j=0}^{m_i-1}z_j^{(i)}.$$
Considering the sequence of reflections determined by $\eta$ that have the same reflection obstacles as the ones determined by $\tilde{\xi}^{(i)}$ for all $i$, we can find at least $(p-2)s_iL$ reflections from the same obstacles as in $\tilde{\xi}^{(i)}$ for all $i$ at the begining of the periodic block and at the end. Thus, by Lemma \ref{5}, we have:
\begin{equation}
\parallel y_{s_iL+j}^{(i)}-z_{s_iL+j}^{(i)}\parallel\le \delta^{s_iL}+\delta^{s_iL}=2\delta^{s_iL},\; j=0,...,(p-2)s_iL-1.
\end{equation}
Now, by (4.1) and (4.5), we have:
\begin{eqnarray}
&&\norm[\bigg]{\rho_\phi(\eta)-\frac{1}{\sum_{i=1}^km_i}\sum_{i=1}^km_i\rho_\phi(\tilde{\xi}^{(i)})}\nonumber\\
&=&\norm[\bigg]{\frac{1}{mpL}\sum_{i=1}^k\sum_{j=0}^{m_i-1}z_j^{(i)}-\frac{1}{mpL}\sum_{i=1}^km_i\left(\frac{1}{m_i}\sum_{j=0}^{m_i-1}y_j^{(i)}\right)}\nonumber\\
&\le&\frac{1}{mpL}\sum_{i=1}^k\sum_{j=0}^{s_iL-1}\parallel y_j^{(i)}-z_j^{(i)}\parallel
+\frac{1}{mpL}\sum_{i=1}^k\sum_{j=s_iL}^{(p-1)s_iL-1}\parallel y_j^{(i)}-z_j^{(i)}\parallel\nonumber\\
&&+\frac{1}{mpL}\sum_{i=1}^k\sum_{j=(p-1)s_iL}^{ps_iL-1}\parallel y_j^{(i)}-z_j^{(i)}\parallel\nonumber\\
&\le&\frac{1}{mpL}\sum_{i=1}^ks_iL2R+\frac{1}{mpL}\sum_{i=1}^k(p-2)s_iL2\delta^{s_iL}
+\frac{1}{mpL}\sum_{i=1}^ks_iL2R\nonumber\\
&\le&\frac{1}{p}(4R+(p-2)2\delta^L\nonumber\\
&\le&\frac{4R}{p}+2\delta^L<\frac{\epsilon}{3}.
\end{eqnarray}
Also, by (4.3), we get:
\begin{eqnarray}
\norm[\bigg]{\frac{1}{mpL}\sum_{i=1}^km_i\rho_\phi(\tilde{\xi}^{(i)})-\sum_{i=1}^kt_i\rho_\phi(\tilde{\xi}^{(i)})}&\le&\sum_{i=1}^k\parallel\rho_\phi(\tilde{\xi}^{(i)})\parallel\norm[\bigg]{\frac{m_i}{mpL}-t_i}\nonumber\\
&\le&kR\frac{\epsilon}{3kR}=\frac{\epsilon}{3}.
\end{eqnarray}
Thus, by (4.4), (4.6) and (4.7), we have:
\begin{eqnarray}
&&\parallel\rho_\phi(\eta)-v\parallel\nonumber\\
&\le&\norm[\bigg]{\rho_\phi(\eta)-\frac{1}{mpL}\sum_{i=1}^km_i\rho_\phi(\tilde{\xi}^{(i)})}+\norm[\bigg]{\frac{1}{mpL}\sum_{i=1}^km_i\rho_\phi(\tilde{\xi}^{(i)})-\sum_{i=1}^kt_i\rho_\phi(\tilde{\xi}^{(i)})}\nonumber\\
&&+\norm[\bigg]{\sum_{i=1}^kt_i\rho_\phi(\tilde{\xi}^{(i)}-\sum_{i=1}^kt_i\rho_i}\nonumber\\
&<&\frac{\epsilon}{3}+\frac{\epsilon}{3}+\sum_{i=1}^kt_i\norm{\rho_\phi(\tilde{\xi}^{(i)})-\rho_i}\nonumber\\
&<&\frac{\epsilon}{3}+\frac{\epsilon}{3}+\sum_{i=1}^kt_i\frac{\epsilon}{3}=\epsilon\nonumber.
\end{eqnarray}
This completes the proof. \endofproof

Now, we will prove that the set $P_\phi$ is dense in $J_\phi$.
\begin{lemma}
\label{4}
If, for some $\xi\in\Sigma$, there exits $\rho_\phi(\xi)=\lim_{n\rightarrow\infty}\frac{1}{n}\sum_{i=0}^{n-1}\phi(\sigma^i\xi)\in\mathbb{R}^N$, then for every $\epsilon>0$ there exists a periodic $\eta\in\Sigma$ such that $\parallel\rho_\phi(\xi)-\rho_\phi(\eta)\parallel<\epsilon$.
\end{lemma}
\emph{Proof}. As $\rho_\phi(\xi)=\lim_{n\rightarrow\infty}\frac{1}{n}\sum_{i=0}^{n-1}\phi(\sigma^i\xi)$, we can find $n_0$ (large) such that 
\begin{equation}\norm[\bigg]{\rho_\phi(\xi)-\frac{1}{n}\sum_{i=0}^{n-1}\phi(\sigma^i\xi)}<\frac{\epsilon}{3},\quad\forall\, n\geq n_0.
\end{equation}
Now take $n>n_0$ and assume that $n=p\ell$, where $p,\ell$ are as in (4.1). Construct $\eta$ to be periodic with period block as follow:
$$\eta=(....;\overbrace{\xi_0,\xi_1,.....,\xi_{n-1}}^{\mbox{periodic block}},....)$$
if $\xi_0\ne\xi_{n-1}$, $$\eta=(....,\overbrace{\xi_0,\xi_1,.....,\xi_n}^{\mbox{periodic block}},....)$$ if $\xi_{n-1}=\xi_0$, then $\xi_n\neq\xi_0$.

We will consider the first case -- the second case is similar. Then $\eta$ is periodic with period $n$.

Let $x_i$ be the reflection points of the billiard trajectory determined by $\xi$, and let $y_i$ be the reflection points of the billiard trajectory determined by $\eta$. Then, by Lemma \ref{1}, $\rho_\phi(\eta)=\frac{1}{n}\sum_{i=0}^{n-1}y_i$. As we can find in $\xi$ and $\eta$ at least $\ell$ reflections from the same obstacles at the begining and the end of the periodic block of $\eta$ corresponding to those of $\xi$, by Lemma \ref{5}, we have:
\begin{equation}
\parallel x_{\ell+j}-y_{\ell+j}\parallel\le2\delta^\ell,\;j=0,...,(p-2)\ell-1.
\end{equation}
Hence, by (4.1) and (4.10), we have:
\begin{eqnarray}
\parallel\rho_\phi(\xi)-\rho_\phi(\eta)\parallel
&=&\norm[\bigg]{\rho_\phi(\xi)-\frac{1}{n}\sum_{i=0}^{n-1}y_i}\nonumber
\end{eqnarray}
\begin{eqnarray}
&\le&\norm[\bigg]{\rho_\phi(\xi)-\frac{1}{n}\sum_{i=0}^{n-1}\phi(\sigma^i\xi)}+\frac{1}{n}\sum_{i=0}^{\ell-1}\parallel x_i-y_i\parallel\nonumber\\
&&+	\frac{1}{n}\sum_{i=\ell}^{(p-1)\ell-1}\parallel x_i-y_i\parallel+\frac{1}{n}\sum_{i=(p-1)\ell}^{p\ell}\parallel x_i-y_i\parallel\nonumber\\
&\le&\frac{\epsilon}{3}+\frac{1}{n}[2\ell R+(p-2)2\ell\delta^\ell+2\ell R]\nonumber\\
&\le&\frac{\epsilon}{3}+\frac{4R}{p}+2\delta^\ell<\epsilon,\nonumber 
\end{eqnarray}
which proves the statement. \endofproof 

\setlength{\parskip}{10pt plus 1pt minus 1pt}

We can now prove our main result.

\medskip
\emph{Proof of Theorem \ref{11}}. The second part follows from Lemma \ref{4}. It remains to prove that the closure $\overline{J_\phi}$ of $J_\phi$ is convex. 

As $\overline{J_\phi}$ is closed, any point in $\overline{J_\phi}$ can be approximated by points in $J_\phi$, so it is enough to prove that $conv(J_\phi)\subset\overline{J_\phi}$. 

Let $\rho_1,\rho_2,...\rho_k$ be elements of $J_\phi$ and let $v=\sum_{i=1}^kt_i\rho_i$ where $t_i>0$ and $\sum_{i=1}^kt_i=1$ (i.e $v$ is a convex combination of $\rho_i,\; i=1,...,k$). Let $\epsilon>0$. For any $i=1,2,...,k$, we can find, by Lemma \ref{4}, a periodic $\eta_i$ such that \begin{equation}\parallel\rho_i-\rho_\phi(\eta_i)\parallel<\frac{\epsilon}{2}.\end{equation} Set $\rho_i'=\rho_\phi(\eta_i),\quad v'=\sum_{i=1}^kt_i\rho_i'.$ By Lemma \ref{3}, we can find a periodic $\tilde{\xi}$ such that \begin{equation}\parallel\rho_\phi(\tilde{\xi})-v\parallel<\frac{\epsilon}{2}.\end{equation}
By (4.11) and (4.12) we have:
\begin{eqnarray}
\norm{\rho_\phi(\tilde{\xi})-v}&=&\norm{\rho_\phi(\tilde{\xi})-v'+v'-v}\nonumber\\
&\le&\norm{\rho_\phi(\tilde{\xi})-v'}+\norm{v'-v}\nonumber\\
&\le&\frac{\epsilon}{2}+\norm{\sum_{i=1}^kt_i\rho_i'-\sum_{i=1}^kt_i\rho_i}\nonumber\\
&\le&\frac{\epsilon}{2}+\sum_{i=1}^k t_i\norm{\rho_i'-\rho_i}\nonumber\\
&\le&\frac{\epsilon}{2}+\sum_{i=1}^k t_i\frac{\epsilon}{2}=\epsilon\nonumber
\end{eqnarray}
Hence $v\in \overline{J_\phi}$, and $\overline{J_\phi}$ is convex. \endofproof

\section{Constructing a trajectory with no rotation vector}
Here, we give an example of a billiard trajectory defined by some $\xi$ in $\Sigma$ such that $\rho_\phi(\xi)$ does not exist.
\setcounter{equation}{0}
\subsection{Example}
Assume that $a_1,a_2,a_3$ are distinct points in $\mathbb{R}^N$, and fix $\epsilon >0$ so small such that $\parallel a_2-a_3\parallel>12\epsilon$, $\parallel a_1-a_2\parallel>2\epsilon$, $\parallel a_1-a_3\parallel>2\epsilon$, $dist(a_1,L)>2\epsilon$, where $L$ is the line $a_2, a_3$. Let $K_1,K_2,K_3$ be strictly convex obstacles containing $a_1,a_2,a_3$, respectively, and having diameter $\le\epsilon$. Then, for any point $x$ in $\partial K_i$, we have 
\begin{equation}\parallel x-a_i\parallel\le\epsilon \end{equation}
Define $\xi\in\Sigma$ as follows:
$$\xi=(.....;B_1,B_2,B_3,.....).$$
where $$B_k=({\underbrace{1,2,....,1,2}_{\mbox{$2^{2k}$}}},{\underbrace{1,3,.....,1,3}_{\mbox{$2^{2k+1}$}}})$$
and $m_k=2^{2k}+2^{2k+1}$ is the length of $B_k$.

We will construct distinct limit points of the sequence
$$\{b_n:b_n=\frac{1}{n}\sum_{j=0}^{n-1}\phi(\sigma^j\xi)\}.$$

Let $n_k=m_1+m_2+.....+m_{k-1}+2^{2k}$ and $\hat{n}_k=m_1+m_2+....+m_k$.

Consider the subsequences $\{b_{n_k}\}$ and $\{b_{\hat{n}_k}\}$ of $b_n$, i.e.
$$\{b_{n_k}:b_{n_k}=\frac{1}{n_k}\sum_{j=0}^{n_k-1}\phi(\sigma^j\xi)\},$$
$$\{b_{\hat{n}_k}:b_{\hat{n}_k}=\frac{1}{\hat{n}_k}\sum_{j=0}^{\hat{n}_k-1}\phi(\sigma^j\xi)\}.$$
Now, for a given $k$, let
$$S_1=\mid\{j:0\leq j\leq n_k,x_j\in\partial K_1\}\mid,$$
$$S_2=\mid\{j:0\leq j\leq n_k,x_j\in\partial K_2\}\mid,$$
$$S_3=\mid\{j:0\leq j\leq n_k,x_j\in\partial K_3\}\mid.$$
Set $$P_i=\, \mid S_i\mid,\quad i=1,2,3.$$
Obviously, $$P_1+P_2+P_3=n_k.$$
From the definitions of $\xi$ and $n_k$, it is clear that half of the reflection points belong to $\partial K_1$ (i.e $P_1=\frac{n_k}{2}$), and so \begin{equation}\norm[\bigg]{\left(\frac{P_1}{n_k}-\frac{1}{2}\right)a_1}=0.\end{equation}
And the other half of $n_k$ will be devided between $\partial K_2$ and $\partial K_3$.

Now 
\begin{eqnarray}
n_k&=&{\underbrace{2^2+2^3}_{\mbox{$m_1$}}}+{\underbrace{2^4+2^5}_{\mbox{$m_2$}}}+.....+{\underbrace{2^{2(k-1)}+2^{2(k-1)+1}}_{\mbox{$m_{2(k-1)}$}}}+2^{2k}\nonumber\\
&=&2^2\left(\frac{2^{2k-1}-1}{2-1}\right)=2^2({2^{2k-1}-1})\nonumber
\end{eqnarray}
From the definitions of $B_k$, $n_k$ and $m_k$, we have that:
\begin{eqnarray}
P_2&=&2+2^3+2^5+.....+2^{2k-1}\nonumber\\
&=&2\left(\frac{2^{2k}-1}{2^2-1}\right)=\frac{2}{3}(2^{2k}-1),\nonumber\end{eqnarray}

and
\begin{eqnarray}
P_3&=&2^2+2^4+....+2^{2k-2}\nonumber\\
&=&2^2\left(\frac{2^{2k-2}-1}{2^2-1}\right)=\frac{4}{3}(2^{2k-2}-1).\nonumber\end{eqnarray}

And so as $k\rightarrow\infty$, we have \begin{equation}\norm[\bigg]{\left(\frac{P_2}{n_k}-\frac{1}{3}\right)a_2}\rightarrow 0, \end{equation}
and
\begin{equation}\norm[\bigg]{\left(\frac{P_3}{n_k}-\frac{1}{6}\right)a_3}\rightarrow 0 \end{equation}
Hence, from (5.1),(5.2),(5.3) and (5.4), we have:
\begin{eqnarray}
\norm[\bigg]{ b_{n_k}-\left(\frac{a_1}{2}+\frac{a_2}{3}+\frac{a_3}{6}\right)}&\leq&\norm[\bigg]{\frac{1}{n_k}\sum_{j\in S_1}(x_j-a_1)+\frac{1}{n_k}\sum_{j\in S_2}(x_j-a_2)+\frac{1}{n_k}\sum_{j\in S_3}(x_j-a_3)}\nonumber\\
&&+\norm[\bigg]{\left(\frac{P_1}{n_k}-\frac{1}{2}\right)a_1}+\norm[\bigg]{\left(\frac{P_2}{n_k}-\frac{1}{3}\right)a_2}+\norm[\bigg]{\left(\frac{P_3}{n_k}-\frac{1}{6}\right)a_3}\nonumber\\
&\le&\frac{P_1}{n_k}\epsilon+\frac{P_2}{n_k}\epsilon+\frac{P_3}{n_k}\epsilon=\epsilon\nonumber
\end{eqnarray}
when k is large enough.

Thus, if $b$ is a limit point of $\{b_{n_k}\}$ (i.e the limit of a convergent subsequence of $\{b_{n_k}\}$), then $$\norm[\bigg]{ b-\left(\frac{a_1}{2}+\frac{a_2}{3}+\frac{a_3}{6}\right)}\le\epsilon.$$ Similarly, it can be proved that if $\hat{b}$ is a limit point of $\{b_{\hat{n}_k}\}$ (i.e the limit of a convergent subsequence of $\{b_{\hat{n}_k}\}$), then $$\norm[\bigg]{ \hat{b}-\left(\frac{a_1}{2}+\frac{a_2}{6}+\frac{a_3}{3}\right)}\le\epsilon.$$

By our assumption, $\parallel a_2-a_3\parallel>12\epsilon$, we have that 
$$\norm[\bigg]{\left(\frac{a_1}{2}+\frac{a_2}{3}+\frac{a_3}{6}\right)-\left(\frac{a_1}{2}+\frac{a_2}{6}+\frac{a_3}{3}\right)}=\norm[\bigg]{\frac{a_2}{6}-\frac{a_3}{6}}>2\epsilon.$$
Thus $b\neq\hat{b}$, so $\lim_{n\rightarrow\infty}b_n$ does not exist. \endofproof

\footnotesize
{\bf Acknowledgements.} Thanks are due to Luchezar Stoyanov for suggesting the problem and continuous support.

\bigskip

\noindent
\textit{School of Mathematics and Statistics\\ University of Western Australia, Perth WA 6009\\ Australia}

zainab.alsheekhhussain@research.uwa.edu.au


\normalsize
\begin{thebibliography}{100}
\bibliographystyle{plain}
\bibitem{1} A. Blokh, Functional rotation numbers for one dimensional maps, Trans. Amer. Math. Soc. 347 (1995), 499-513.
\bibitem{12} A. Blokh, M. Misiurewicz and N. Simanyi, Rotation sets of billiards with one obstacle, Commun. Math. Phys. 266 (2006), 239-265. 
\bibitem{2} W. Geller and M. Misiurewicz, Rotation and entropy, Trans. Amer. Math. Soc. 351 (1999), 2927-2948.
\bibitem{3} M. Ikawa, Decay of solutions of the wave equation in exterior of several strictly convex bodies, Ann. Inst. Four. 4 (1988), 113-146.
\bibitem{28} O. Jenkinson, Rotation, entropy, and equlibrium maps, Trans. Amer. Math. Soc. 353 (2001), 3713-3739.
\bibitem{5} T. Kucherenko and C. Wolf, Geometry and entropy of generalized rotation sets, (2012), arXiv:1210.0135vI[math.Ds].
\bibitem{30}D. Lind and B. Marcus, An Introduction to Symbolic Dynamics and Coding, Cambridge University Press, 1995. 
\bibitem{24} M. Misiurewicz, K. Ziemian, Rotation sets and ergodic measures for torus homeomorphisms, Fund. Math. 137 (1991), 45-52.
\bibitem{25} M. Misiurewicz, K. Ziemian, Rotation sets for maps of tori, J. London Math. Soc. 40 (1989), 490-509.
\bibitem{6} M. Misiurewicz, Rotation theory, Scholarpedia 2.10 (2007):3873.
\bibitem{7} T. Morita, The symbolic representation of billiards without boundary condition, Trans. Amer. Math. Soc. 325 (1991), 819-828.
\bibitem{8} V. Petkov and L. Stoyanov, Geometry of Reflecting Rays and Inverse Spectral Problems, Wiley, Chichester (1992).
\bibitem{9} Ya. Sinai, Dynamical systems with elastic reflections: ergodic properties of dispering billiards, Uspehi Mat. Nauk 10 (1970), 141-192 (in Russian); English transl.: Russian Math. Surv. 25 (1970), 137-189.
\bibitem{10*}J. Sj\"{o}strond, Geometric bounds on the density of resonances for semiclassical problems, Duke Math. J. 60 (1990), 1-57.
\bibitem{20} L. Stoyanov, An estimate from above of the number of periodic orbits for semi-dispersed billiards, Commun. Math. Phys. 124 (1989), 217-227.

\bibitem{21} L. Stoyanov, Exponential instability for a class of dispering billiards, Ergod. Th. \& Dynam. Sys. 19 (1999), 201-226.
\bibitem{29} P. Walters, An Introduction to Ergodic Theory, Springer, Berlin (1982).
\bibitem{11} K. Ziemian, Rotation sets for subshifts of finite type, Fund. Math. 146 (1995), 189-201.
\end{thebibliography}
\end{document}